\documentclass[a4paper,12pt]{amsart}
\usepackage{amsfonts,amsmath,amssymb,amscd,latexsym}

\newtheorem{thm}{Theorem}[section]

\newtheorem{lem}[thm]{Lemma}
\newtheorem{cor}[thm]{Corollary}

\newcommand{\la}{\langle\,}
\newcommand{\ra}{\,\rangle}
\newcommand{\be}[1]{\begin{equation}\label{#1}}
\newcommand{\ee}{\end{equation}}
\newcommand{\Aut}{\mathop{\mathrm{Aut}}}
\newcommand{\zz}{{\mathbb Z}}
\newcommand{\qq}{{\mathbb Q}}
\newcommand{\ve}{\varepsilon}
\def\A{\mathcal A}
\def\R{\mathcal R}
\def\L{\mathcal L}
\def\M{\mathcal M}
\def\N{\mathcal N}
\def\F{\mathcal F}
\def\ninf{\N_{\mathrm{inf}}}
\def\nfin{\N_{\mathrm{fin}}}

\begin{document}

\hfuzz=80pt
{

\title{Growth of positive words and lower bounds of the growth rate
for Thompson's groups $F(p)$}

}

\author{Jos\'e Burillo}
\address
{Escola Polit\`{e}cnica Superior de Castelldefels, UPC, Avda del
Canal Ol\'{\i}mpic s/n, 08860 Castelldefels, Barcelona, Spain}
\email
{burillo{@}mat.upc.es}

\author{Victor Guba}
\address
{Vologda State Pedagogical University,
6 S.\,Orlov Street,
Vologda
160600
Russia}
\email
{guba{@}uni-vologda.ac.ru}

\thanks{Both authors thank the Centre de Recerca Matem\`{a}tica
for the hospitality and support. The research of the second author
was also partially supported by the RFFI grant 05-01-00895.}

\date{\today}

%\dedicatory{}

\begin{abstract}
Let $F(p)$, $p\ge2$ be the family of generalized Thompson's groups.
Here $F(2)$ is the famous Richard Thompson's group usually denoted
by $F$. We find the growth rate of the monoid of positive words
in $F(p)$ and show that it does not exceed $p+1/2$. Also we describe
new normal forms for elements of $F(p)$ and, using these forms, we
find a lower bound for the growth rate of $F(p)$ in its natural
generators. This lower bound asymptotically equals $(p-1/2)\log_2 e+1/2$
for large values of $p$.
\end{abstract}

\maketitle

\section*{Introduction}

The family of generalized Thompson's groups $F(p)$ was introduced
by K. S. Brown in \cite{Br87}. Additional facts about these groups
can be found in \cite{BCS,Stein}. The case $p=2$ corresponds to the
famous Richard Thompson's group $F$. See the survey \cite{CFP} for
details about this group.

The groups $F(p)$ have many common features. All of them are embeddable
into each other \cite{BrGuz}. None of them has free non-abelian
subgroups. None of these groups satisfy any nontrivial group law.
The derived subgroups of each of the $F(p)$ is simple (infinitely
generated). Every proper homomorphic image of $F(p)$ is abelian
(so these groups are not residually finite). Each $F(p)$ is finitely
presented and has quadratic Dehn function \cite{Gu05}.

Each of these groups has a faithful representation by piecewise linear
functions. The word problem has an easy solution in each of these
groups. Also all these groups are diagram groups in the sense of
\cite{GuSa}. Namely, $F(p)$ is a diagram group over a very simple
semigroup presentation $\la x\mid x=x^p\ra$. It follows then from
\cite[Section 15]{GuSa} that $F(p)$ has solvable conjugacy problem.
Each group $F(p)$ satisfies homological finiteness condition
$\F_\infty$. All integer homology groups $H_n(F(p),\zz)$ are free
abelian of finite rank and the Poincar\'e series are rational
\cite{GS}.

However, there is some difference between the groups of this family.
Brin \cite{Brin} described the group $\Aut F$ for $F=F(2)$. Some
information about automorphisms of $F(p)$, where $p>2$, can be
found in \cite{BrGuz}, where it is shown that already for $p=3$
there are ``wild" automorphisms of $F(p)$.
\vspace{1ex}

The goal of this article is to obtain analogs of some results
for the group $F$. The first author found the growth function
of the monoid of positive elements of $F$. This function is
rational, namely, it equals
$$
\frac{1-x^2}{1-2x-x^2+x^3}.
$$

\noindent
Notice that the elements $x_0$, $x_1$, \dots, $x_{p-1}$ generate a free
submonoid of rank $p$ in $F(p)$. Thus the growth rate of positive
elements in $F(p)$ is at least $p$. In this paper we show that for
any $p$, the exact value of the growth rate of positive elements is
only slightly higher than $p$ --- it never exceeds $p+1/2$.

Guba and Sapir \cite{GuSa97} found two new normal forms for elements of
$F$. One of them is a normal form in the infinite set of generators.
This normal form is locally testable (unlike the standard normal form).
It has the same feature as the normal form in the free group: a word is
in a normal form if and only if all its subwords of length $2$ are in
the normal form. In this paper, we find such a form for every $F(p)$.
Another normal form constructed in \cite{GuSa97} for $F$ allows one to
construct a regular set of normal forms in $F$. We find an analogous
construction for each $F(p)$.

Using the above regular normal form, the second author proved in
\cite{Gu04} that the growth rate of the group $F$ in generators
$x_0$, $x_1$ is at least $(3+\sqrt{5})/2$. Notice that neither the
growth function, nor the growth rate for $F$ is known at the present.
In this paper we find a lower bound of the growth rate for each
of the groups $F(p)$, where the generating set consists of $x_0$,
$x_1$, \dots, $x_{p-1}$. We show that the lower bound is a root of a
certain algebraic equation and find the asymptotic behaviour of this
root. For large values of $p$, this is $(p-1/2)\log_2 e+1/2$, where
$\log_2 e=1.442695\ldots$\,.

The plan of the paper is as follows. In Section \ref{prel} we recall
the definition of the family $F(p)$ of generalized Thompson's groups
and some basic facts about growth functions and growth rates. This
Section also contains a description of (positive) elements in $F(p)$
in terms of rooted $p$-trees.

In Section \ref{Ford} we describe Fordham's method to calculate the word
length in $F(p)$. We restrict ourselves to the case of positive words
only (the description for this case is much simpler). Recall that for
the case $p=2$, a fast algoritm to find the word length metric was
described in \cite{Fordham,Fordham2}. This algorithm is very effective
but it has quite a complicated description. A simplification of the
method due to Belk and Brown can be found in \cite{BeBr}. One of the
easiest algorithms to find the word length in $F$ (the so-called Length
Formula) is contained in \cite[Section 5]{Gu04}. Notice that for $p>2$,
none of the simplified versions exists so we use Fordham's approach
from \cite{FordP}.

In Section \ref{gfgs}, using Fordham's method, we find equations for
generating functions describing the growth of $F_+(p)$. We solve
these equations in Section \ref{gfgrfp+} and show that the generating
function for positive words in $F(p)$ is irrational provided $p>2$
(unlike the case $p=2$). Then we find the growth rate of positive
words in $F(p)$ as a root of an algebraic equation. We prove that
this growth rate never exceeds $p+1/2$ approaching this value as
$p$ approaches infinity. Thus the set $F_+(p)$ of all positive words
is not much higher than the free submonoid generated by $x_0$, $x_1$,
\dots, $x_{p-1}$.

Section \ref{nnfefp} describes two new normal forms of elements in
$F(p)$. The first of these forms is locally testable (one needs to test
only subwords of length $2$, similar to a free group). The second of
the normal forms leads to a regular language that represents each
element of $F(p)$ exactly once. Based on that regular language, we
construct the corresponding automaton and find a lower bound for the
growth rate of $F(p)$ in Section \ref{lbgrfp}. This lower bound is
given as a root of an algebraic equation. We also describe its
asymptotic behaviour showing that it approaches $(p-1/2)\log_2 e+1/2$
for large values of $p$.

\section{Preliminaries}
\label{prel}

The family of generalized Thompson's group can be defined as follows.
The group $F(p)$ is the group of all piecewise linear self homeomorphisms
of the unit interval $[0,1]$ that are orientation preserving (that is,
send $0$ to zero and $1$ to $1$) with all slopes integer powers of $p$
and such that their singularities (breakpoints of the derivative) belong
to $\zz[\,\frac1p\,]$. The group $F(p)$ admits a presentation given by
\be{presfp}
\la x_i\ (i\ge0)\,\mid x_jx_i=x_ix_{j+p-1}\ (i<j)\ra.
\ee

\noindent
This presentation is infinite, but a close examination shows that the
group is actually finitely generated, since $x_0$, $x_1$, \dots, $x_{p-1}$
are sufficient to generate it. In fact, the group is finitely presented.
The finite presentation is awkward and it is not used much. The symmetric
and simple nature of the infinite presentation makes it much more adequate
for almost all purposes.

One such example where the infinite presentation is particularly
appropriate is in the construction of the normal form. A word given in the
generators $x_i$ and their inverses, can have its generators moved around
according to the relators, and the result is the following well-known
statement:

\begin{thm}
\label{stnf}
An element in $F(p)$ always admits an expression of the form
$$
x_{i_1}x_{i_2}\cdots x_{i_m}x_{j_n}^{-1}\cdots x_{j_2}^{-1}x_{j_1}^{-1},
$$
where
$$
i_1\le i_2\le\cdots\le i_m,\ j_1\le j_2\le\cdots\le j_n.
$$
\end{thm}

In general, this expression is not unique, but for every element
there is a unique word of this type which satisfies certain
technical condition (see \cite{CFP} for details). This unique word
is called the {\em standard normal form\/} for the element of
$F(p)$.

Observe that the infinite presentation for $F(p)$ is actually a monoid
presentation. Hence $F(p)$ admits a submonoid, the submonoid $F_+(p)$
given by the same presentation, whose elements are called
{\em positive words\/}. Theorem \ref{stnf} shows that $F(p)$ is the
group of right fractions of this monoid.

An element of $F(p)$ can be represented by two subdivisions of the
interval $[0,1]$, namely, the subdivision into intervals which get mapped
linearly to each other. A subdivision of this type, where the dividing
points are all in $\zz[\,\frac1p\,]$, can always be obtained by subsequent
subdivisions of the interval into $p$ equal pieces. Hence, a subdivision
of the interval is equivalent to a rooted tree where each vertex has
valence $p+1$ except the root, which has valence $p$ (or $1$ in case when
the tree consists of the root only), and the {\em leaves\/}, which have
valence $1$. A node (except the root and the leaves) is pictured to have
one edge going up and $p$ edges going down to its $p$ {\em children\/}.
These trees will be called {\em rooted\/} $p$-{\em trees\/}. An element
of $F(p)$ is then represented by a pair of rooted $p$-trees called the
{\em source\/} tree and the {\em target\/} tree. This representation has
been extensively studied in the case $p=2$. Note that positive words can
be represented by a single $p$-tree, because the other tree is always the
same: the tree which consists of all right carets.

A piece of these $p$-trees consisting of a node and its $p$ edges going
down to its children is called a {\em caret\/}. Carets are the building
blocks of the trees and they give rise to the algorithm for finding the
word metric in $F(p)$, see Section \ref{Ford}.
\vspace{2ex}

As stated in the introduction, the exact growth function for the
groups $F(p)$ is not known. In this paper we will give lower bounds
for growth rates of these groups, computing lower bounds for the
number of elements in each length.

To be precise, given a finitely generated group $G$ with finite
generating set $X$, denote its sphere of radius $n$ by
$$
{\mathbf S}(n)=\{\,g\in G\mid \ell(g)=n\,\},
$$
where $\ell(g)$ is the length of $g\in G$ in the set of generators $X$.
We also have the ball of radius $n$
$$
{\mathbf B(n)}=\bigcup_{k=0}^n{\mathbf S}(k).
$$
If $\gamma_n=\#{\mathbf B}(n)$, the series
$$
\Gamma(x)=\sum_{n=0}^\infty\gamma_nx^n
$$
is called the (general) {\em growth function\/} for $G$ with respect to
$X$, and the number
$$
\gamma=\lim_{n\to\infty}\gamma_n
$$
is the {\em growth rate\/} of $G$ with respect to $X$. The limit
always exists due to the submultiplicative property of $\gamma_n$,
that is, $\gamma_{m+n}\le\gamma_m\gamma_n$ for all $m,n\ge0$. Also,
the {\em spherical growth function\/} is given by
$\sigma_n=\#{\mathbf S}(n)$ and
$$
\Sigma(x)=\sum_{n=0}^\infty\sigma_nx^n,
$$
which has the same growth rate as the general growth function (for all
infinite groups). For details about growth functions, see, for instance,
\cite{GriHa}.

If $P\subseteq G$ is a subset of a group, not necessarily a subgroup,
we can define the growth functions of the set $P$ by the same formulas
as above but where the coefficients are actually the cardinals of the
sets $P\cap{\mathbf B}(n)$ or $P\cap{\mathbf S}(n)$. The goal for one
of the next sections is to compute the growth series of the subset
$F_+(p)$ in $F(p)$. In order to do that, we need to describe the
algorithm for calculating the word metric in $F(p)$.

\section{Positive words in Thompson's groups $F(p)$ and Fordham's method}
\label{Ford}

In 1995, S. Blake Fordham \cite{Fordham} constructed an algorithm which,
for any given element in $F=F(2)$, finds its distance to the identity in
the word metric given by generators $x_0$, $x_1$. This algorithm consists
in defining different types of carets, then having each caret of the
source tree paired to its corresponding caret in the target tree, and
assigning a weight to each type of pairs of carets. A table is given for
all possible pairs of types, with the assignment of the weight. The sum
of all the weights of all the pairs is the exact distance from the element
to the identity. In a set of unpublished notes \cite{FordP}, Fordham
extends his method to the groups $F(p)$. This method will be the starting
block of the computation.

The method used to compute this growth will be an extension to $F(p)$ of
the method developed in \cite{Burillo} for the case of $F=F(2)$. Consider
a positive element of $F(p)$. As we know, the element can be represented
by a rooted $p$-tree. We are going to define different types of carets
and their weights, following Fordham \cite{FordP}.
\vspace{1ex}

A caret will be called {\em left\/} or {\em right\/} if it is
situated in the leftmost edge of the tree or in the rightmost edge,
and {\em middle\/} or {\em interior\/} if it is situated in the
middle, i.e. if it is not right or left. For instance, a caret is
left if it represents a subinterval of $[0,1]$ which has left
endpoint equal to zero. Middle carets will be subdivided into $p-1$
types, denoted by $\M^1$, $\M^2$, \dots, $\M^{p-1}$ according to
which caret they are children of, and its position as child.

The children of a caret are subdivided in two types, the
predecessors and the successors. This subdivision will give a total
order to the set of carets, with a caret being always after its
predecessor children and before its successors. The definitions of
the caret types are as follows:

\begin{itemize}
\item The root caret is special. Its children are:
\begin{itemize}
\item Its left child is a left caret and it is the only predecessor.
\item Its middle children are successors, and have types
$\M^1$, $\M^2$, \dots, $\M^{p-2}$, in order-preserving way.
\item Its right child is obviously a successor and a right caret.
\end{itemize}
\item A left caret has the following children:
\begin{itemize}
\item Its only predecessor is the left child, a left caret.
\item All the other children are successors, all middle carets, and of
types $\M^1$, $\M^2$, \dots, $\M^{p-1}$, in order.
\end{itemize}
\item A right caret has the following children:
\begin{itemize}
\item One single predecessor of type $\M^{p-1}$.
\item It has $p-1$ successors, which in order are of types $\M^1$, $\M^2$,
\dots, $\M^{p-2}$ and the last one of type $\R$.
\end{itemize}
\item A caret of type $\M^i$ ($1\le i\le p-1$) has the following
children:
\begin{itemize}
\item The first $p-i$ children are predecessors, and their types are
$\M^i$, \dots, $\M^{p-1}$.
\item The other $i$ children are successors, and they are of types
$\M^1$, $\M^2$, \dots, $\M^i$.
\end{itemize}
\end{itemize}

For the purposes of computing the length of an element, these caret types
are subdivided in further types depending on the existence of predecessor
and successor types. This classification is actually more complicated in
Fordham's paper but we do not need the total strength of the method since
we are dealing only with positive words. We will indicate also which is
the weight of each caret for the purposes of the computation of the length
of a positive word.

The caret types are as follows:

\begin{itemize}
\item The root, which has always weight zero.
\item Left carets, which have always weight one.
\item Carets of type $\R_{\varnothing}$ are right carets whose all
successors are right carets, i.e., it has no middle successors. Its
only successors hang from its rightmost leaf. These carets carry
weight zero.
\item Carets of type $\R_M$ are right carets which are not
$\R_{\varnothing}$, that is, which have middle successors. Observe that
the middle successors do not have to be immediate successors, they can
be successors of successors. Carets of type $\R_M$ have weight two.
\item Carets of type $\M_{\varnothing}^i$ are middle carets which do not
have any successor children. They carry weight one.
\item Carets of type $\M_M^i$ are middle carets which have at least a
successor child. These carets have weight three.
\end{itemize}

%
%PICTURE OF CARET OF TYPE M^i GOES HERE (maybe)
%

Observe that the index on the middle carets is only necessary to identify
its successors, but it has no role in the weight assignment beyond that
one.

Now, the main theorem giving the length is as follows:

\begin{thm}
{\rm (S. B. Fordham) \cite{FordP}} Given a positive word in $F(p)$
represented by a rooted $p$-tree, the distance from this element to
the identity $($in the word metric for $F(p)$ with generators $x_0$,
$x_1$, \dots, $x_{p-1})$ is equal to the total sum of the weights of
its carets.
\end{thm}

\section{Generating functions for the growth of positive words}
\label{gfgs}

Once the theorem for the length has been established, now the computation
of the growth function is reduced to a combinatorial problem, namely,
finding how many trees have a given weight, according to the rules above.
The method for finding the number of trees with a given weight is to
split the trees in several ones in such a way that recurrences can be
found. The reader can see details about generating functions in
\cite{Wilf}, and can see this method used already in \cite{Burillo}.

We will make use of several sequences:

\begin{itemize}
\item The sequence $s_n=\#(F_+(p)\cap{\mathbf S}(n))$. This is the number
of trees which have weight $n$.
\item The sequence $l_n$. This sequence gives the number of subtrees
which can be left subtrees of a rooted $p$-tree and such that its total
weight is $n$. The subtrees are required to be strict, that is, the main
tree does not qualify as a left subtree.
\item Analogously the sequence $r_n$ is the sequence of possible right subtrees
of weight $n$.
\item The sequence $m^{(i)}_n$ for $i=1,\ldots,p-1$, gives the number of
interior subtrees which start with a caret of type $\M^i$. Observe that
this subtree is completely composed of middle carets, and also with
total weight $n$.
\end{itemize}

Observe that the subtrees are always considered as subtrees of the main
tree, which means that, for instance, a left subtree never has carets of
type $\R$ because that would mean it is the total tree. A subtree which
starts in an $\M^i$ caret has all interior carets.

Each one of these sequences will have its generating function:
$$
S(x)=\sum_{n=0}^\infty s_nx^n\qquad
L(x)=\sum_{n=0}^\infty l_nx^n\qquad
R(x)=\sum_{n=0}^\infty r_nx^n\qquad
M_i(x)=\sum_{n=0}^\infty m_n^{(i)}x^n.
$$

Now we will establish relations between the sequences which will give
functional equations for their generating functions, which then will
allow us to find the growth of the submonoid of positive words. For
instance, if one considers the tree representing a word, and assumes the
tree has total weight $n$, since the root has weight zero, the weight has
to be distributed among all the $p$ children subtrees. Hence, a tree of
total weight $n$ will be obtained every time that we take a family of
subtrees such that the sum of their separate weights as subtrees is $n$.

This fact gives the first formula satisfied by the sequences, and also by
the generating functions:
\be{sn}
s_n=\sum_{j_0+\cdots+j_{p-1}=n}l_{j_0}m^{(1)}_{j_1}
\cdots m^{(p-2)}_{j_{p-2}}r_{j_{p-1}}
\ee
\be{gfs}
\qquad S=LM_1\cdots M_{p-2}R.
\ee

To find a formula for the function $L(x)$ of left subtrees, one needs to
consider that left carets have weight $1$. Hence the different subtrees
only have to add up to $n-1$. The formula is
$$
l_n=\sum_{j_0+\cdots+j_{p-1}=n-1}l_{j_0}m^{(1)}_{j_1}
\cdots m^{(p-2)}_{j_{p-2}}m^{(p-1)}_{j_{p-1}}
$$
\be{gfl}
L-1=xL M_1 M_2\cdots M_{p-1}.
\ee
The formula for the generating functions is obtained by multiplying each
side of the formula for sequences by $x^n$. The right hand side has an
$x$ multiplying because the indices are shifted by one.

For the function for right trees, one has to take into account the
fact that a right caret can be of type $\R_{\varnothing}$ or $\R_M$,
with weights zero and two respectively. For the first possibility,
the caret is of type $\R_{\varnothing}$, and all its successors have
no weight. Observe that in a positive word there can be one and only
one caret of type $\R_{\varnothing}$, because any others would be
reducible. Hence, if the caret is of type $\R_{\varnothing}$, all
the weight is concentrated in its only predecessor. So there are as
many right subtrees of this type as trees of the type $\M^{p-1}$
with the same weight, which gives the first part of the recurrence
equal to $m^{(p-1)}_n$.

If the right caret is of type $\R_M$, it carries weight $2$ and one the
successors is necessarily nonempty with a middle caret somewhere. Hence
if one of the successors is necessarily nonempty, the term in the
recurrence has all possible weights for these successors. The formula is
$$
r_n=m^{(p-1)}_n+\sum_{\substack{j_0+\cdots+j_{p-1}=n-2\\
j_1+\cdots+j_{p-1}\ge1}}m^{(p-1)}_{j_0}m^{(1)}_{j_1}\cdots
m^{(p-2)}_{j_{p-2}}r_{j_{p-1}}
$$
\be{gfr}
R=M_{p-1}+x^2(M_1M_2\cdots M_{p-1}R-M_{p-1}).
\ee

Finally, the middle subtrees are the ones whose children are also middle
subtrees and hence facilitate the resolution of the equations. A middle
caret of type $\M^i$ has either weight $1$ if its successors are empty or
weight $3$ if one of the successor subtrees is nonempty. Both cases
correspond to the two adding terms of the formula for the sequence:
$$
m^{(i)}_n=\sum_{j_i+\cdots+j_{p-1}=n-1} m^{(i)}_{j_i}\cdots
m^{(p-1)}_{j_{p-1}}+\sum_{\substack{j_0+j_1+\cdots+j_{p-1}=n-3\\
j_{p-i}+\cdots+j_{p-1}\ge1}}m^{(i)}_{j_0}\cdots
m^{(p-1)}_{j_{p-i-1}}m^{(1)}_{j_{p-i}}\cdots m^{(i)}_{j_{p-1}}
$$
which gives the following formula for the generating functions:
\be{gfm}
M_i-1=xM_iM_{i+1}\cdots M_{p-1}+x^3M_iM_{i+1}
\cdots M_{p-1}(M_1\cdots M_{i-1}M_i-1).
\ee
Solving these equations will give us information on the function $S(x)$,
which is the one we are interested in, and the growth of positive
elements in the groups $F(p)$.

\section{Growth functions and growth rates of $F_+(p)$}
\label{gfgrfp+}

Now we collect formulas (\ref{gfs}), (\ref{gfl}), (\ref{gfr}),
(\ref{gfm}) to find the equation on $S(x)$ and the radius of
convergence of the corresponding series. First of all, we have to
mention that $F_+(p)$ has a free submonoid generated by $x_0$,
$x_1$, \dots, $x_{p-1}$ and so the growth rate of $F_+(p)$ is at
least $p$. As we will see at the end of this Section, the exact
value of the growth rate is only slightly larger than $p$. (In fact,
it is always less than $p+1/2$.)

Let
$$
M(x)=M_1(x)M_2(x)\cdots M_{p-1}(x).
$$

\begin{lem}
\label{mmm} For all $0\le i\le p-1$, we have
$$
M_1M_2\cdots M_i=\frac{x^{-2}}{(1-x^3M)^i}+1-x^2.
$$
\end{lem}

\proof
We proceed by induction on $i$. If $i=0$, then the result
is obvious. Let $1\le i\le p-1$. Formula (\ref{gfm}) can
be written as
$$
M_i=1+\frac{xM}{M_1\cdots M_{i-1}}+
x^3\left(M_iM-\frac{M}{M_1\cdots M_{i-1}}\right).
$$
Therefore,
$$
M_1\cdots M_i=M_1\cdots M_{i-1}+xM+x^3M\cdot M_1\cdots M_i-x^3M
$$
and so
$$
M_1\cdots M_i(1-x^3M)=(x-x^3)M+M_1\cdots M_{i-1}.
$$
Using the inductive assumption, we have
$$
M_1\cdots M_i(1-x^3M)=(x-x^3)M+\frac{x^{-2}}{(1-x^3M)^{i-1}}+1-x^2=
\frac{x^{-2}}{(1-x^3M)^{i-1}}+(1-x^2)(1-x^3M).
$$
Now the only thing left to do is to divide by $1-x^3M$.
\endproof

Taking $i=p-1$ gives us

\begin{cor}
\label{eqonm}
The function $M=M(x)$ satisfies
$$
x^2M=\frac1{(1-x^3M)^{p-1}}+x^2-1.
$$
\end{cor}

Now we express $S(x)$ in terms of $M(x)$. It follows from
(\ref{gfl}) and (\ref{gfr}) that
$$
L=\frac1{1-xM}\qquad R=\frac{(1-x^2)M_{p-1}}{1-x^2M}.
$$
Now, using (\ref{gfs}), we have
\be{S}
S=\frac{LMR}{M_{p-1}}=\frac{(1-x^2)M}{(1-xM)(1-x^2M)}.
\ee

The first author proved in \cite{Burillo} that the growth
function $S(x)$ of positive elements of $F=F(2)$ is rational
(although $M(x)$ is irrational). Now we have the following

\begin{thm}
\label{irr}
The growth function $S(x)$ of positive elements in $F(p)$
is irrational provided $p\ge3$.
\end{thm}

\proof
Let $N=(1-x^3M)^{-1}$. From Corollary \ref{eqonm} we have
$$
1-N^{-1}=x^3M=xN^{p-1}+x^3-x.
$$
Hence $N=N(x)$ satisfies the equation
\be{eqonn}
xN^p+(x^3-x-1)N+1=0.
\ee

Suppose that $S(x)$ is rational. Then it follows from (\ref{S}) that
$M(x)$ satisfies a quadratic equation with coefficients in the field
$\qq(x)$ of rational functions. Since $M=x^{-3}(1-N^{-1})$, the
function $N(x)$ also satisfies an equation of degree at most $2$
over $\qq(x)$. This implies that the polynomial
$f(t)=xt^p+(x^3-x-1)t+1$ from $\qq(x)[t]$ is divisible by a
polynomial of degree at most $2$. Since $p\ge3$, the polynomial
$f(t)$ is reducible over $\qq(x)$. A standard algebraic trick (using
Gauss' lemma) implies that $f(t)$ is a product of two polynomials
from $\zz[x][t]$ of degree less than $p$. Taking $x=1$, we obtain
that the polynomial $t^p-t+1$ is reducible over $\qq$. However, this
contradicts a result from \cite{Lj}.
\endproof

Now we will find the growth rate of $F_+(p)$. To do that, we need to
take the radius of convergence of the series for $S(x)$ and take the
reciprocal. Observe that from (\ref{sn}) we deduce $m_n\le s_n$ for
all $n\ge0$. This implies that
$$
(\limsup_{n\to\infty}m_n)^{-1}\ge(\limsup_{n\to\infty}s_n)^{-1},
$$
that is, the radius of convergence of the series $S(x)$ does not
exceed the one for the series $M(x)$. Let $x>0$ be a real
number such that $S(x)$ converges. Then $M(x)$ also converges
and formula (\ref{S}) holds.

To find the radius of convergence of $S(x)$, we need to find
the smallest positive real number such that the denominator
of the right hand side of (\ref{S}) is zero. Since $M(x)$ is
increasing and $0<x<1$, the smallest positive solution of the
equation $M(x)=x^{-1}$ will not exceed the smallest positive
solution of the equation $M(x)=x^{-2}$. Therefore, we need to
solve the equation $M(x)=x^{-1}$. Notice that $M(x)$ increases
and $x^{-1}$ decreases so we can just speak about a positive
root of this equation. Using (\ref{eqonm}), we get
$x=(1-x^2)^{-(p-1)}+x^2-1$, that is, we need to find the
positive root of
\be{xx}
(1-x^2)^{p-1}(1+x-x^2)=1.
\ee
The growth rate of $F_+(p)$ will thus be equal to $x^{-1}$. We
already know that the growth rate of $F_+(p)$ is at least $p$,
as it was mentioned in the beginning of this Section. Hence
$x\le1/p$.

Let us rewrite this equation in the following form:
$$
p-1=\frac{\ln(1+x-x^2)}{-\ln(1-x^2)}.
$$
From the Taylor formula for $\ln (1+y)$, we deduce the inequality
$$
y-y^2/2<\ln(1+y)<y-y^2/2+y^3/3,
$$
where $y>0$, and then we get
$\ln(1+x-x^2)<x-3x^2/2+4x^3/3-3x^4/2+x^5-x^6/3<x-3x^2/2+4x^3/3$.
Since $-\ln(1-x^2)>x^2$, we have $p-1<x^{-1}-3/2+4x/3\le
x^{-1}-3/2+4/3p$. So $x^{-1}>p+1/2-4/3p=p+1/2+o(1)$ as $p\to\infty$.

Now we want to show that $x^{-1}<p+1/2$. We have
$\ln(1+x-x^2)>x-x^2-(x-x^2)^2/2=x-3x^2/2+x^3-x^4/2$ and
$-\ln(1-x^2)=x^2+x^4/2+x^6/3+\cdots<x^2+x^4(1+x^2+x^4+\cdots)/2=
x^2+x^4/(2-2x^2)\le x^2+2x^4/3$ because $x\le1/p\le1/2$. This gives
$p-1>(x-3x^2/2+x^3-x^4/2)/(x^2+2x^4/3)=(1-3x/2+x^2-x^3/2)/(x+2x^3/3)$.
Finally,
$$
p-1>\frac{1-3x/2+x^2-x^3/2}{x+2x^3/3}=\frac1x-\frac{9-2x+3x^3}{2(3+2x^2)}>
1/x-3/2
$$
since $3x^2-6x-2<0$ on $[0;1]$. This gives $x^{-1}<p+1/2$, as
desired. So we get the following result.

\begin{thm}
\label{grratef+}
The growth rate of the monoid $F_+(p)$ of positive elements
in the group $F(p)$ generated by $x_0$, $x_1$, \dots, $x_{p-1}$
is a number $\zeta_p$, which is the root of equation
$$
(y^2-1)^{p-1}(y^2+y-1)=y^{2p}.
$$
This number has the form $\zeta_p=p+\lambda_p$, where
$0<\lambda_p<1/2$ for all $p$ and $\lambda_p\to 1/2$ as
$p\to\infty$.
\end{thm}

Indeed, we proved inequalities $p+1/2-4/3p<x^{-1}<p+1/2$,
where $x$ is the solution of (\ref{xx}). The inequality
$x^{-1}>p$ obviously follows for $p\ge3$; if $p=2$, then
it is known from \cite{Burillo} that $\zeta_2>2.24$.

The equation in the statement of Theorem \ref{grratef+} is
equivalent to (\ref{xx}) via the substitution $y=1/x$. Notice that
$x$ and $y$ are roots of polynomials of degree $2p-1$ with integer
coefficients. Also let us mention without proof that $\lambda_p$ is
strictly increasing with respect to $p$.

The number $\zeta_p$ gives a lower bound for the growth rate of the
group $F(p)$. However, this estimate can be essentially improved.

\section{New normal forms for elements of $F(p)$}
\label{nnfefp}

We are going to find two new normal forms for elements of $F(p)$.
They will be analogs of the normal forms constructed in \cite{GuSa97}
for the case $F=F(2)$.

The first of these normal forms will involve the infinite set of
generators $\Sigma=\{\,x_i\ (i\ge0)\,\}$. Consider the following
rewriting system $\Gamma=\Gamma(p)$ over the alphabet
$\Sigma^{\pm1}=\Sigma\cup\Sigma^{-1}$ (basic facts about rewriting
systems can be found in \cite{BO,DJ}):

\begin{enumerate}
\item
$x_i^\ve x_i^{-\ve}\to1$\quad($i\ge0$, $\ve=\pm1$)
\item
$x_j^\ve x_i\to x_ix_{j+p-1}^\ve$\quad($j>i$, $\ve=\pm1$)
\item
$x_{j+p-1}^\ve x_i^{-1}\to x_i^{-1}x_j^\ve$\quad($j>i$, $\ve=\pm1$)
\end{enumerate}

Notice that for every rewriting rule of $\Gamma$, the left hand
side and the right hand side are equal in $F(p)$.

It is easy to see that $\Gamma$ is {\em terminating\/}, that is,
for every word $w$, the process of applying rewriting rules to $w$
always terminates. Indeed, $\Gamma$ either decreases the length
of a word or it preserves the length. In the second case, if
we make a vector that consists of subscripts of a word, the
rewriting rules will decrease this vector lexicographically.

Since $\Gamma$ is terminating, applying the rewriting rules to a
word $w$ gives us a word $v$ that cannot be reduced (that is, no
more rewriting rules can be applied to $v$). We say that $v$ is an
{\em irreducible form\/} of $w$. Now we are going to check that
$\Gamma$ is also {\em confluent\/}, that is, every word has a unique
irreducible form. To do that, we apply the Diamond Lemma. In our
case, this means that if we have rewriting rules of the form $ab\to
u$, $bc\to v$, where $a$, $b$, $c$, $d$ are letters and $u$, $v$ are
words, then $uc$ and $av$ have a common descendant. There are only
finitely many cases to check, and all of them are easy. We will show
one of these cases, the rest is left to the reader.

Let us take the rewriting rules
$x_{k+p-1}^\ve x_j^{-1}\to x_j^{-1}x_k^\ve$ and
$x_j^{-1}x_i\to x_ix_{j+p-1}^{-1}$, where $k>j>i$, $\ve=\pm1$.
We have:
$$
x_j^{-1}x_k^\ve x_i\to x_j^{-1}x_ix_{k+p-1}^\ve\to
x_ix_{j+p-1}^{-1}x_{k+p-1}^\ve
$$
and
$$
x_{k+p-1}^\ve x_ix_{j+p-1}^{-1}\to x_ix_{k+2p-2}^\ve x_{j-p-1}^{-1}\to
x_ix_{j+p-1}^{-1}x_{k+p-1}^\ve.
$$
So the words have a common descendant.

Now we know that $\Gamma$ is {\em complete\/}, that is, terminating
and confluent. Therefore, each element of $F(p)$ can be uniquely
represented by an irreducible word. So we have proved the following

\begin{thm}
\label{nfinf}
Each element $g\in F(p)$ can be uniquely represented as a word
of the form
$$
N(g)=x_{i_1}^{\ve_1}x_{i_2}^{\ve_2}\cdots x_{i_m}^{\ve_m},
$$
where $m\ge0$, $\ve_1,\ve_2,\ldots,\ve_m=\pm1$, and for every $1\le k<m$
one of the following conditions holds:
\begin{itemize}
\item
$i_k<i_{k+1}$
\item
$i_k=i_{k+1}$ and $\ve_k=\ve_{k+1}$
\item
$0<i_k-i_{k+1}<p$ and $\ve_{k+1}=-1$.
\end{itemize}
\end{thm}

Indeed, the conditions listed in the statement exactly mean that
the word $N(g)$ is irreducible, that is, it has no subwords that
are left hand sides of the rewriting rules of $\Gamma$. The set
of these irreducible words over $\Sigma^{\pm1}$ will be denoted
by $\ninf$.

Notice that the set $\ninf$ has the following property: a word
belongs to $\ninf$ if and only if all its subwords of length $2$
belong to $\ninf$. That is, the normal form of Theorem \ref{nfinf}
is locally testable.
\vspace{1ex}

Now we will construct another normal form for elements of $F(p)$.
Now all words will involve only the finite set of generators
$x_0^{\pm1}$, $x_1^{\pm1}$, \dots, $x_{p-1}^{\pm1}$. Moreover, these
normal forms will give a regular language closed under taking
subwords. Notice that this gives a regular spanning tree in the
Cayley graph of $F(p)$ in the above generators. As in \cite{GuSa97}
for the case $p=2$, this tree is not geodesic.

It is possible to write down a new rewriting system in order to get
the normal form we wish to construct. However, it will take too much
effort to prove that the rewriting system ijs complete. We choose an
approach that differs from \cite{GuSa97}.

Let $j\ge1$. Then $j$ can be uniquely expressed in the form
$j=r+d(p-1)$, where $1\le r\le p-1$, $d\ge0$. In this case $x_j$
equals in $F(p)$ to the word $x_0^{-d}x_rx_0^d$. For any word $w$
over $\Sigma^{\pm1}$, replace each letter of the form $x_j^\ve$
($j\ge1$, $\ve=\pm1$) by $x_0^{-d}x_rx_0^d$, where $j=r+d(p-1)$,
$1\le r\le p-1$, $d\ge0$ and then freely reduce all subwords of the
form $x_0^\ve x_0^{-\ve}$ ($\ve=\pm1$). We obtain a word in
generators $x_0^{\pm1}$, $x_1^{\pm1}$, \dots, $x_{p-1}^{\pm1}$
denoted by $\bar w$.

\begin{lem}
\label{forb}
If $w\in\ninf$, then $\bar w$ has no subwords of the following
form:
\begin{enumerate}
\item
$x_i^\ve x_i^{-\ve}$\quad$(0\le i\le r-1)$
\item
$x_\alpha^\ve x_0^kx_\beta$\quad$(k\ge0$, $1\le\beta<\alpha\le r-1)$
\item
$x_\alpha^\ve x_0^{k+1}x_\beta^{-1}$\quad
$(k\ge0$, $1\le\beta<\alpha\le r-1)$
\item
$x_\alpha^\ve x_0^{k+1}x_\beta$\quad$(k\ge0$, $1\le\alpha\le\beta\le r-1)$
\item
$x_\alpha^\ve x_0^{k+2}x_\beta^{-1}$\quad
$(k\ge0$, $1\le\alpha\le\beta\le r-1)$
\end{enumerate}
\end{lem}

The words of the form 1) -- 5) are called {\em forbidden subwords\/}.
The set of words in $\{\,x_0^{\pm1},x_1^{\pm1},\ldots,x_{p-1}^{\pm1}\,\}$
without forbidden subwords will be denoted by $\nfin$.

\proof
Let $w\in\ninf$ have the form
\be{wnf}
w=x_0^{k_0}x_{\alpha_1}^{l_1}x_0^{k_1}x_{\alpha_2}^{l_2}\cdots
x_0^{k_{h-1}}x_{\alpha_h}^{l_h}x_0^{k_h},
\ee
where $h\ge0$, $\alpha_i=r_i+d_i(p-1)$, $1\le r_i\le p-1$, $d_i\ge0$,
$l_i\ne0$ for all $1\le i\le h$. By definition,
\be{barw}
\bar w=x_0^{k_0-d_1}x_{r_1}^{l_1}x_0^{d_1+k_1-d_2}x_{r_2}^{l_2}\cdots
x_0^{d_{h-1}+k_{h-1}-d_h}x_{r_h}^{l_h}x_0^{d_h+k_h}.
\ee
Suppose that $\bar w$ is not freely irreducible. Then there
exist an $i$ from $1$ to $h-1$ such that $d_i+k_i-d_{i+1}=0$,
$r_i=r_{i+1}$, and and $l_il_{i+1}<0$. By definition, words
from $\ninf$ have no subwords of the form $x_j^{\pm1}x_0$ for
$j>0$ and also have no subwords of the form $x_j^{\pm1}x_0^{-1}$
for $j\ge p$. This implies $k_i\le0$.

Suppose that $k_i<0$. Then $d_i=0$ (otherwise $\alpha_i\ge p$).
Since $d_{i+1}\ge0$, we obtain $d_i+k_i-d_{i+1}=k_i-d_{i+1}<0$.
This is a contradiction. Therefore, $k_i=0$ and so $d_i=d_{i+1}$.
This implies $\alpha_i=r_i+d_i(p-1)=r_{i+1}+d_{i+1}(p-1)=\alpha_{i+1}$.
Thus the word $w$ is not freely irreducible since $l_il_{i+1}<0$. We
have a contradiction. This proves that $\bar w$ has no subwords of
the form 1).

Suppose that $\bar w$ has a subword of one of the forms 2) -- 5).
Let $x_{r_i}^{\pm1}x_0^{d_i+k_i-d_{i+1}}x_{r_{i+1}}^{\pm1}$ be
such a subword, where $1\le i<h$. As above, $k_i\le0$. Suppose that
$k_i\ne0$. This implies $d_i=0$ and $d_i+k_i-d_{i+1}<0$. But none
of the words 2) -- 5) can contain $x_0^{-1}$. This allows us to
conclude that $k_i=0$ and $\bar w$ contains
$v=x_{r_i}^{\pm1}x_0^{d_i-d_{i+1}}x_{r_{i+1}}^{\pm1}$ as a subword.

Suppose that $v$ satisfies condition 2). This means that $d_i\ge d_{i+1}$,
$r_i>r_{i+1}$, $l_{i+1}>0$. Hence $w$ contains
$x_{\alpha_i}^{\pm1}x_{\alpha_{i+1}}$, where
$\alpha_i=r_i+d_i(p-1)>r_{i+1}+d_{i+1}(p-1)$. So $w$ does not belong
to $\ninf$, which is impossible.

Suppose that $v$ satisfies condition 3). Now $d_i-d_{i+1}\ge1$,
$r_i>r_{i+1}$, $l_{i+1}<0$. This leads to
$\alpha_i-\alpha_{i+1}=(r_i-r_{i+1})+(d_i-d_{i+1})(p-1)\ge p$,
which also contradicts $w\in\ninf$.

Suppose that $v$ satisfies condition 4). Then $d_i-d_{i+1}\ge1$,
$r_i\le r_{i+1}$, $l_{i+1}>0$. Now $r_i-r_{i+1}\ge1-(p-1)=2-p$
and so
$\alpha_i-\alpha_{i+1}=(r_i-r_{i+1})+(d_i-d_{i+1})(p-1)\ge (p-1)+(2-p)>0$.
Thus $w$ contains $x_{\alpha_i}^{\pm1}x_{\alpha_{i+1}}$ with
$\alpha_i>\alpha_{i+1}$. This cannot happen by definition of $\ninf$.

Finally, suppose that $v$ satisfies condition 5). Now we have
$d_i-d_{i+1}\ge2$ and so
$\alpha_i-\alpha_{i+1}=(r_i-r_{i+1})+(d_i-d_{i+1})(p-1)\ge2(p-1)+(2-p)=p$.
However, it should be $\alpha_i-\alpha_{i+1}<p$ because $w\in\ninf$.

The proof is complete.
\endproof

For every $g\in F(p)$, we have the word $\overline{N(g)}\in\nfin$
that represents $g$. We will prove that $g$ is represented uniquely
by a word from $\nfin$. This will follow from

\begin{lem}
\label{uniq}
The mapping $w\mapsto\bar w$ from $\ninf$ to $\nfin$ is a bijection.
\end{lem}

\proof We prove first that the mapping $w\mapsto\bar w$ from $\ninf$
to $\nfin$ is injective. As above, let $w\in\ninf$ have the form
(\ref{wnf}). Thus $\bar w$ equals (\ref{barw}). Suppose that we know
the word $\bar w$, that is, we know the numbers $m_0=k_0-d_1$,
$m_1=d_1+k_1-d_2$, \dots, $m_{h-1}=d_{h-1}+k_{h-1}-d_h$,
$m_h=k_h+d_h$. Our aim is to recover the numbers $k_0$, $d_1$,
$k_1$, \dots, $d_{h-1}$, $k_{h-1}$, $d_h$, $k_h$.

Let $h\ge1$. It follows from the definition of $\ninf$ that
$k_h\le0$. Moreover, either $k_h<0$ and $d_h=0$, or $k_h=0$.
In the first case $m_h=k_h+d_h=k_h<0$, in the second case
$m_h=k_h+d_h=d_h\ge0$. Since we know $m_h$, we can distinguish
between these two cases. Namely, if $m_h<0$, then $d_h=0$, $k_h=m_h$.
If $m_h\ge0$, then $k_h=0$, $d_h=m_h$. Now we know $d_h$ and $k_h$.

If $h\ge2$, then $k_{h-1}\le0$. As above, we have one of the two
cases: $k_{h-1}<0$, $d_{h-1}=0$, or $k_{h-1}=0$. The number
$k_{h-1}+d_{h-1}$ is negative in the first case and nonnegative
in the second case. But this number equals $m_{h-1}+d_h$, so we
know it and thus we are able to distinguish these cases. In the
first case we have $d_{h-1}=0$, $k_{h-1}=m_{h-1}+d_h$; in the
second case --- $k_{h-1}=0$, $d_{h-1}=m_{h-1}+d_h$. Therefore,
we know $d_{h-1}$ and $k_{h-1}$.

Continuing in this way, we get the values of $d_{h-2}$, $k_{h-2}$,
\dots, $d_1$, $k_1$. At the final step we get $k_0=m_0+d_1$.
\vspace{1ex}

Now we show that the mapping is surjective. We start with a word
from $\nfin$. This word has the form \be{surj}
x_0^{m_0}x_{r_1}^{l_1}x_0^{m_1}x_{r_2}^{l_2}\cdots
x_0^{m_{h-1}}x_{r_h}^{l_h}x_0^{m_h}. \ee Using the rules described
in the first part of the proof, we define the numbers $k_0$, $d_1$,
$k_1$, \dots, $d_h$, $k_h$. It follows that $d_i\ge0$ for all $i$
from $1$ to $h$. So we can form a word $w$ as in (\ref{wnf}), where
$\alpha_i=r_i+d_i(p-1)$ ($1\le i\le h$). It is obvious that $\bar w$
equals the word (\ref{surj}). It remains to prove that $w$ belongs
to $\ninf$.

Let us assume the contrary. Since $w$ has no subwords of the
form $x_0^\ve x_0^{-\ve}$, it should contain one of the
following subwords:

a) $x_i^{\ve}x_i^{-\ve}$\quad($i\ge1$, $\ve=\pm1$);

b) $x_j^{\pm1}x_i$\quad($j>i$);

c) $x_{j+p-1}^{\pm1}x_i^{-1}$\quad($j>i$).

In case a), $\bar w$ will contain a forbidden subword of the form
$x_r^\ve x_r^{-\ve}$. Notice that $k_i\le0$ for all $1\le i\le h$;
if $k_i<0$, then $d_i=0$. This means that in cases b) and c) one
has $i\ge1$. Let $j=\alpha+d(p-1)$, $i=\beta+d'(p-1)$, where
$1\le\alpha,\beta\le r-1$, $d,d'\ge0$. Applying the ``bar" mapping
to b) and c), we see that the word $\bar w$ contains
$u=x_\alpha^{\pm1}x_0^{d-d'}x_\beta$ in case b) and
$v=x_\alpha^{\pm1}x_0^{d-d'+1}x_\beta^{-1}$ in case c). Since $j>i$,
we have $(d-d')(p-1)>\beta-\alpha>-(p-1)$. Hence $d-d'\ge0$. If $u$
is not forbidden, then $\alpha\le\beta$. But in this case $d-d'>0$
so $u$ has to be forbidden anyway. If $v$ is not forbidden, then
we also have $\alpha\le\beta$, which implies $d-d'+1\ge2$. We
have a final contradiction.

The proof is complete.
\endproof

From Lemmas \ref{forb} and \ref{uniq} we obtain

\begin{thm}
\label{nffin}
Each element $g\in F(p)$ can be uniquely represented by a word
$w\in\nfin$. This means that for every $g\in F(p)$ there is exactly
one word over $\{\,x_0^{\pm1},x_1^{\pm1},\ldots,x_{p-1}^{-1}\,\}$
that represents $g$ and has no forbidden subwords. This gives a
regular set of normal forms for the group $F(p)$.
\end{thm}

Indeed, the set of forbidden subwords is a regular language. So the
set $\nfin$ of words that do not contain forbidden subwords will be
also regular. Throughout the rest of the paper, we will denote this
language by $\L_p$. (For basic properties of regular languages see
\cite{Sal}.)

\section{Lower bounds for the growth rates of $F(p)$}
\label{lbgrfp}

A lower bound of $(3+\sqrt{5})/2=2.618\ldots$ for the growth rate of
$F=F(2)$ was obtained by the second author in \cite{Gu04}. Now we
will find a similar lower bound for each $F(p)$. In the previous
section, we constructed a regular language $\L_p$ of normal forms
for $F(p)$. Each word of length $n$ in $\L_p$ is at a distance at
most $n$ from the identity in the Cayley graph of $F(p)$. So the
growth function of $\L_p$ does not exceed the number of elements in
the ball of radius $n$ for $F(p)$. Then, finding the growth function
and the growth rate of $\L_p$, we find a lower bound for the growth
rate of the group $F(p)$.

An automaton to recognize the language $\L_p$ has $3p+1$ states.
However, it is easier to construct a directed graph with only $2p+1$
vertices (states). This graph will be denoted by $\A_p$ and we will
also call it an automaton although its edges have no labels. The
description of $\A_p$ is as follows.

The vertices (states) of $\A_p$ are denoted by $q$, $q_0$, $q_1$,
\dots, $q_{p-1}$, $q_{1,0}$, $q_{2,0}$, \dots, $q_{p-1,0}$, $\bar
q$. They will correspond to the following partition of $\L_p$ into
disjoint subsets:
\begin{itemize}
\item
The set $\{\,1\,\}$ that consists of the empty word (state $q$).
\item
The set of words that end with $x_0^{\pm1}$ and do not have a
terminal segment of the form $x_i^{\pm1}x_0^k$, where $1\le i\le p-1$,
$k\ge1$ (state $q_0$).
\item
The set of words that end with $x_i^{\pm1}$ (state $q_i$ for
each $1\le i\le p-1$).
\item
The set of words that end with $x_i^{\pm1}x_0$ (state $q_{i,0}$ for
each $1\le i\le p-1$).
\item
The set of words that end with $x_i^{\pm1}x_0^k$ for some $1\le i\le p-1$
and $k\ge2$ (state $\bar q$).
\end{itemize}

Let $w\in\L_p$. If $w$ is empty, then $wx_i^{\pm1}$ will be in
$\L_p$ for all $0\le i\le p-1$. We draw two arrows from $q$ to
$q_i$ for each $0\le i\le p-1$.

Let $w$ correspond to the state $q_0$. Then $w=vx_0^\ve$ for
some word $v$ and for some $\ve=\pm1$. The word $wx_0^{\ve}$
will be in $\L_p$; for each $1\le i\le p-1$ the word $wx_i^{\pm1}$
will be also in $\L_p$ since $w$ has no terminal segments of the
form $x_i^{\pm1}x_0^k$ ($1\le i\le p-1$, $k\ge1$. Thus we draw an
arrow from $q_0$ to itself and two arrows from $q_0$ to each
$q_i$ ($1\le i\le p-1$).

Let $w$ correspond to $q_i$ ($1\le i\le p-1$). The words $wx_0^{-1}$
and $wx_0$ belong to $\L_p$; we draw an arrow from $q_i$ to $q_0$
and an arrow from $q_i$ to $q_{i,0}$. The words $wx_j$ belong to
$\L_p$ whenever $i<j\le p-1$; the words $wx_j^{-1}$ belong to $\L_p$
for all $1\le j\le p-1$. So we draw one arrow from $q_i$ to each
$q_1$, \dots, $q_i$ and two arrows from $q_i$ to each $q_{i+1}$,
\dots, $q_{p-1}$.

Let $w$ correspond to $q_{i,0}$ ($1\le i\le p-1$). The word
$wx_0$ is in $\L_p$ and we draw an arrow from $q_{i,0}$ to $\bar q$.
Also $wx_j^{-1}\in\L_p$ whenever $i\le j\le p-1$. So one arrow
goes from $q_{i,0}$ to each $q_i$, \dots, $q_{p-1}$. No other
arrows can appear.

Finally, let $w$ correspond to $\bar q$. Now only $wx_0$ leads to a
word in $\L_p$; it corresponds to an arrow from $\bar q$ to itself.

The description of $\A_p$ is complete. Notice that the number of
words in $\L_p$ of length $n$ is exactly the number of (directed)
paths of length $n$ in $\A_p$ starting at $q$. We would like to
compute the number of paths in $\A_p$ of length $n$ starting at $q$
and ending at a given state. For each state we consider the
corresponding generating function. Namely, to each vertex $v$ we
assign a series of the form $\sum_{n=0}^\infty a_nt^n$, where $a_n$
is the number of paths in $\A_p$ starting at $q$ and ending at $v$.
These generating functions will be denoted by $f$, $f_i$ ($0\le i\le
p-1$), $f_{i,0}$ ($1\le i\le p-1$), $\bar f$ for each of the states,
respectively. We will write down a system of equations for these
functions.

First of all, it is clear that $f(t)=1$. To find $f_0$, we mention
that two arrows go from $q$ into $q_0$ and one arrow goes into $q_0$
from each of the states $q_0$, $q_1$, \dots, $q_{p-1}$. Hence
\be{fq0}
f_0=t(2f+f_0+f_1+\cdots+f_{p-1}).
\ee

Given a vertex $q_i$ ($1\le i\le p-1$), we observe that two arrows
go into $q_i$ from $q$, $q_0$, \dots, $q_{i-1}$ and one arrow from
$q_i$, \dots, $q_{p-1}$. Also one arrow goes into $q_i$ from each
$q_{1,0}$, \dots, $q_{i,0}$. Thus
\be{fqi}
f_i=
t(2f+2f_0+\cdots+2f_{i-1}+f_i+\cdots+f_{p-1})+t(f_{1,0}+\cdots+f_{i,0}).
\ee

Notice that $f_{i,0}=tf_i$ for each $1\le i\le p-1$ because only
one arrow goes into $q_{i,0}$ (from the state $q_i$). Thus we can
rewrite (\ref{fqi}) as follows:
\be{fqi+}
f_i=
t(2f+2f_0+\cdots+2f_{i-1}+f_i+\cdots+f_{p-1})+t^2(f_1+\cdots+f_i).
\ee

Finally, there is one arrow that goes into $\bar q$ from each of
the states $q_{1,0}$, \dots, $q_{p-1,0}$, $\bar q$. So
\be{fbarq}
\bar f=t(f_{1,0}+\cdots+f_{p-1,0}+\bar f)=t^2(f_1+\cdots+f_{p-1})+t\bar f.
\ee

In order to solve the system, let us consider the difference of
equation (\ref{fqi+}) with $i=1$ and (\ref{fq0}). This gives
$f_1-f_0=tf_0+t^2f_1$, that is, $(1-t^2)f_1=(1+t)f_0$. So
$f_0=(1-t)f_1$.

Now suppose that $1\le i<p-1$. If we take the difference between
$f_{i+1}$ and $f_i$ using (\ref{fqi+}), we obtain
$f_{i+1}-f_i=tf_i+t^2f_{i+1}$, which implies $f_i=(1-t)f_{i+1}$.

Now for all $0\le i\le p-1$ one has $f_i=(1-t)^{p-1-i}f_{p-1}$. The
equation (\ref{fq0}) becomes
$$
(1-t)^{p-1}f_{p-1}=2t+tf_{p-1}((1-t)^{p-1}+\cdots+(1-t)+1)=
2t+f_{p-1}(1-(1-t)^p).
$$
This gives
\be{frc}
f_{p-1}(t)=\frac{2t}{(1-t)^p+(1-t)^{p-1}-1}.
\ee

In order to find the number of words in $\L_p$ having length $n$,
we need to add all the generating functions for all states. The
result will be
$$
\Phi_p(t)=1+f_0+f_1+\cdots+f_{p-1}+t(f_1+\cdots+f_{p-1})+
\frac{t^2}{1-t}\cdot(f_1+\cdots+f_{p-1})
$$
(here we used (\ref{fbarq}) to express $\bar f$). Taking
into account that
$$
f_1+\cdots+f_{p-1}=\frac{1-t}t\cdot f_0-2
$$
from (\ref{fq0}), we have
$$
\Phi_p(t)=1+f_0+\left(1+t+\frac{t^2}{1-t}\right)(f_1+\cdots+f_{p-1})=
1+\left(1+\frac1t\right)f_0-\frac2{1-t}.
$$
Now, using $f_0=(1-t)^{p-1}f_{p-1}$ and (\ref{frc}), we finally have
\be{phi}
\Phi_p(t)=
\frac{1+t}{1-t}\cdot\frac{1-t(1-t)^{p-1}}{(1-t)^p+(1-t)^{p-1}-1}\,.
\ee
This is the generating function for $\L_p$. Thus the growth rate of $\L_p$
will be the reciprocal of $t$, where $t$ is the smallest positive root
of the denominator of the right hand side of (\ref{phi}).

The number $y=(1-t)^{-1}$ is the root of $y^p=y+1$. It is clear that
$y>1$. Let $y=1+x$, where $x>0$. We would like to solve the equation
$(1+x)^p=2+x$. Notice that $0<x<1$. Since $(1+x)^p<3$, the root $x$
approaches $0$ as $p$ goes to infinity.

The equation $(1+x)^p=2+x$ can be written as
$$
p=\frac{\ln(2+x)}{\ln(1+x)}=\frac{\ln 2+\ln(1+x/2)}{\ln(1+x)}=
\frac{\ln 2+x/2+o(x)}{x-x^2/2+o(x^2)}.
$$
Therefore,
\be{xinv}
\frac{p-\frac12}{\ln2}=\frac{1+o(x)}{x-x^2+o{x^2}}=x^{-1}(1+x/2+o(x))=
x^{-1}+\frac12+o(1)
\ee
as $p\to\infty$.

We are interested in the number $\xi_p=t^{-1}$, where $t$ is the
root of $(1-t)^p+(1-t)^{p-1}-1=0$. Here $(1-t)^{-1}=y$, where
$y=1+x$ is the root of $y^p=y+1$. It is easy to see that $\xi_p=1+x^{-1}$.
So we deduce from (\ref{xinv}) that
$$
\xi_p=\frac{p-\frac12}{\ln 2}+\frac12+o(1),\quad p\to\infty.
$$
It is also easy to see that $\xi_p=t^{-1}$ satisfies the following
equation: $(2\xi-1)(\xi-1)^{p-1}=\xi^p$. So we proved

\begin{thm}
\label{lowbound} The growth rate of the group $F(p)$, $p\ge2$ has a
lower bound of $\xi_p$, where $\xi_p$ satisfies the equation
$$
(2\xi-1)(\xi-1)^{p-1}=\xi^p.
$$
The following asymptotic formula holds:
$$
\xi_p=\frac{p-\frac12}{\ln 2}+\frac12+o(1),\quad p\to\infty.
$$
\end{thm}

Here are several numerical values of $\xi_p$:
$$
\xi_2=2.618033989
$$
$$
\xi_3=4.079595623
$$
$$
\xi_4=5.530132718
$$
$$
\xi_5=6.977144180
$$
and so on. For large values of $p$, the growth rate of $F(p)$
is at least $0.72(2p-1)$ (recall that $2p-1$ is the maximum
value of the growth rate of a $p$-generated group; this happens
if and only if the group is free of rank $p$).
\vspace{1ex}

It would be interesting to find nontrivial upper bounds for the
growth rates of $F(p)$. This means to find a constant $c<1$ such
that the growth rate of $F(p)$ in its natural generators does not
exceed $c(2p-1)$.

\end{document}